\scrollmode
\documentclass[12pt]{amsart}
\usepackage{graphicx}
\usepackage{amssymb}
\usepackage{psfrag}
\usepackage{mathrsfs}
\usepackage{subfigure}
\usepackage{cite}

\newif\iffigs
\figstrue
\figsfalse

\newcommand{\sC}{\mathscr{C}}
\newcommand{\sR}{\mathscr{R}}
\newcommand{\sP}{\mathscr{P}}
\newcommand{\sA}{\mathscr{A}}
\newcommand{\sH}{\mathscr{H}}

\newcommand{\Z}{\mathbb{Z}}

\newcommand{\weylposrts}{\Phi^{+}}
\newcommand{\affrt}{\tilde{\alpha}}
\newcommand{\Waff}{\tilde{W}}
\newcommand{\Saff}{\tilde{S}}
\newcommand{\Reduced}{\mathtt{Red}}
\newcommand{\join}[2]{{#1\!\!\relbar\protect\joinrel\protect\relbar\!\!#2}}
\renewcommand{\tilde}[1]{\widetilde{#1}}
\newcommand{\innprod}[2]{(#1, #2)}
\newcommand{\ggt}{>\!\!>}

\renewcommand{\emptyset}{\varnothing}

\DeclareMathOperator{\GL}{GL}

\swapnumbers
\theoremstyle{plain}
\newtheorem{theorem}[subsection]{Theorem}%[section]
\newtheorem{conjecture}[subsection]{Conjecture}%[section]
\newtheorem{proposition}[subsection]{Proposition}%[section]
\newtheorem{lemma}[subsection]{Lemma}%[section]
\newtheorem{corollary}[subsection]{Corollary}%[section]
 
\theoremstyle{definition}
\newtheorem{definition}[subsection]{Definition}%[section]
\theoremstyle{remark}
\newtheorem{remark}[subsection]{Remark}%[section]
\date{15 July 2008}
\title{Automata and cells in affine Weyl groups}
\author{Paul E. Gunnells}
\address{Department of Mathematics and Statistics\\
University of Massachusetts\\
Amherst, MA 01003}
\email{gunnells@math.umass.edu}


\begin{document}
\begin{abstract}
Let $\Waff$ be an affine Weyl group, and let $C$ be a left, right, or
two-sided Kazhdan--Lusztig cell in $\Waff$.  Let $\Reduced (C)$ be the
set of all reduced expressions of elements of $C$, regarded as a
formal language in the sense of the theory of computation.  We show
that $\Reduced (C)$ is a regular language.  Hence the reduced
expressions of the elements in any Kazhdan--Lusztig cell can be
enumerated by a finite state automaton.
\end{abstract}

\maketitle 

%%%%%%%%%%%%%%%%%%
%                %
%  introduction  %
%                %
%%%%%%%%%%%%%%%%%%
\section{Introduction}\label{s:intro}

\subsection{} Let $W$ be a Coxeter group with generating set $S$.  In
their work on Coxeter groups and Hecke algebras, Kazhdan and Lusztig
defined a partition of $W$ into sets called \emph{cells}.  When $W$ is
a Weyl or affine Weyl group, it is known that cells have deep
connections with many areas of algebra and geometry, such as
singularities of Schubert varieties \cite{kl2}, representations of
$p$-adic groups \cite{lusztig83}, characters of finite groups of Lie
type \cite{lusztig84}, and the geometry of unipotent conjugacy classes in
simple complex algebraic groups \cite{bezru1, bezru2}.

The definition of cells is quite complicated (\S\ref{ss:celldefs}).
It involves the construction of a subtle equivalence relation on $W$
built from both easy and difficult combinatorial data.  In particular
from the definition it is not clear how ``computable'' cells are.  For
instance, it is highly nontrivial to decide whether two elements of
$W$ lie in the same cell or not, or to characterize all elements in a
given cell.  Nevertheless, in all known examples where cells have been
explicitly computed, one sees that cells ultimately have a relatively
simple geometric and combinatorial structure.  We refer to
\cite{cass, notices} for examples and further discussion of this
phenomemon. 

\subsection{} This paper addresses the following computational
problem: given a cell $C$ in a Coxeter group $W$, how can we encode
the (typically) infinite amount of data represented by $C$ with a
finite structure?  In this generality, this question was first
considered by Casselman, who phrased an answer in terms of
\emph{finite state automata} (\S\ref{ss:automata}).  More precisely,
let $\Reduced (W)$ be the set of all reduced expressions of all
elements of $W$, considered as a subset of the free monoid on the
generating set $S$.  We regard $\Reduced (W)$ as a formal languange in
the sense of the theory of computation.  Given a cell $C$, let
$\Reduced (C)\subset \Reduced (W)$ be the set of all reduced
expressions of all $w\in C$.  Then we have the following conjecture of
Casselman:

\begin{conjecture}\label{conj:casselman}
For any Coxeter group $W$ and any cell $C\subset W$, there exists a
finite state automaton accepting the language $\Reduced (C)$.  That
is, the language $\Reduced (C)$ is regular.
\end{conjecture}

The precise definitions of the terms in Conjecture
\ref{conj:casselman} are given in \S\ref{ss:automata}; here we give an
informal sense of what the conjecture means.

A finite state automaton is a simple theoretical model of a computer.
It has finite memory and can complete only one task:
acceptance/rejection of its input.  More precisely, given a word
$a_{1}\dotsb a_{k}$ on some alphabet, an automaton reads the word from
left to right, and while doing so moves through finitely many memory
states.  After reading the word, the automaton decides based on which
state it occupies whether or not to accept the word or throw it away.
A language is called regular if one can find a finite state
automaton accepting exactly the words in the language.

Thus Casselman's conjecture implies that given any cell $C$, there
exists a simple machine $\sA (C)$ that decides whether or not $w\in W$
lies in $C$ simply by reading through a reduced expression
$s_{1}\dotsb s_{k}$ for $w$.  The finiteness of the automaton implies
that the decision is only based on finitely many patterns appearing
in the expression.  One can also use $\sA (C)$ to systematically list
all reduced expressions of all elements of $C$.  

We remark that work of Brink--Howlett \cite{bh} shows that the
language $\Reduced (W)$ of all reduced expressions of all elements of
$W$ is regular, although this does not prove Conjecture
\ref{conj:casselman}: a sublanguage of a regular language need not be
regular.  Indeed, it is not even clear how one can use the tools
underlying the fundamental results of Brink--Howlett to attack
Conjecture \ref{conj:casselman}.  Moreover, Conjecture
\ref{conj:casselman} does not specify what structures in $W$ should be
used to build the machines $\sA (C)$.

Casselman's conjecture is trivially true for any finite Coxeter group,
in particular for Weyl groups, since for such groups the language
$\Reduced (W)$ is obviously finite.  The first infinite example of
Conjecture \ref{conj:casselman} follows from work of Shi \cite{shi} on
Kazhdan--Lusztig cells and Eriksson \cite{eriksson} and Headley
\cite{headley} on automata.  More precisely, let $W=\tilde{A}_{n}$,
the affine Weyl group of type $A$.  Shi showed that $W$ can be
partitioned into finitely many geometrically defined subsets, called
\emph{sign-type regions}, such that each left cell $C$ of $W$ is a
union of finitely many such regions.  Headley showed that the
sign-type regions can be used as a set of states for an automaton
$\sA$ recognizing $\Reduced (W)$.  Together these results imply
Conjecture \ref{conj:casselman} for $\tilde{A}_{n}$, since given $C$
one can modify $\sA $ to only accept the reduced expressions
corresponding to elements of $C$.

\subsection{} In this paper we prove Conjecture \ref{conj:casselman}
when the Coxeter group is an affine Weyl group $\Waff$ (Theorem
\ref{thm:main}).  The proof uses two ingredients.

The first is a family of finite state automata $\sA_{N}$, $N\in
\Z_{\geq 0}$, each of which recognizes $\Reduced (\Waff)$.  The
construction generalizes work of Eriksson \cite{eriksson} and Headley
\cite{headley}.  Each $\sA_{N}$ is built from the complement of a
certain affine hyperplane arrangement $\sH_{N}$.

The second is a result of Du \cite{du.polyhedra}, who proved that each
left cell $C$ of $\Waff$ can be represented as the union of a finite
set of convex polyhedra of a certain type.  We show that if $N\ggt 0$,
then we can write each of these polyhedra as a finite union of regions
in the complement of $\sH_{N}$.  This allows us to define $\sA (C)$ by
identifying the set of states of $\sA_{N}$ that correspond exactly to the
reduced expressions of elements of $C$.  

\subsection{Acknowledgements} We thank M.~Belolipetsky, C.~Bonnafe,
W.~Casselman, J.~Guilhot, and J.~Humphreys for helpful comments.  The
results in this paper were discovered through computation using
modified versions of various programs due to W.~Casselman, D.~Holt,
and F.~du~Cloux.  This work was partially supported by the NSF through
grants DMS 04--01525, 06--19492, 08--01214.

%%%%%%%%%%%%%%%%
%              %
%  background  %
%              %
%%%%%%%%%%%%%%%%
\section{Background}\label{s:background}

\subsection{}
In this section we recall background and standard notation.  For more
details and proofs we refer to \cite{humph.book, bb} for Coxeter
groups, \cite{kl, bb} for Kazhdan--Lusztig cells, and \cite{aho} for
automata and formal languages.

Let $\Phi$ be an irreducible, reduced root system, and let $W$ be the
associated Weyl group.  Decompose $\Phi$ into a union of positive and
negative roots $\Phi^{+}\cup \Phi^{-}$.  Let $\Delta \subset \Phi^{+}$
be the simple roots and let $S\subset W$ be the corresponding subset
of generators in the presentation of $W$ as a Coxeter group.  For each
$s\in S$, we write $\alpha_{s}\in \Delta$ for the associated simple
root.  Conversely, given a simple root $\alpha$, we write
$s_{\alpha}\in S$ for the associated generator.

We assume $\Phi$ spans a real vector space $V$ equipped with a
$W$-invariant inner product $\innprod{\phantom{a}}{\phantom{a}}$.  The
group $W$ acts on $V$ as usual: if $\alpha \in \Phi$, then we have the
reflection
\[
v\longmapsto v-\frac{2\innprod{v}{\alpha}}{\innprod{\alpha}{\alpha}}\alpha, 
\]
which maps $W$ faithfully onto a subgroup of $\GL (V)$.  We will write
this reflection action as a \emph{right} action: $v\mapsto v\cdot
s_{\alpha}$.  Let $\sC^{+}$ be the positive Weyl chamber determined
by $\Delta$.

The root system $\Phi$ has a unique highest root $\affrt$.  By
definition $\affrt$ has the property that for any $\beta \in \Phi$,
the difference $\affrt -\beta$ can be written as a nonnegative linear
combination of the simple roots.  One also knows \cite[IV, \S 1.8,
Prop. 25(iv)]{bourbaki} that for any positive root $\beta$, we have
\[
\frac{2\innprod{\beta}{\affrt}}{\innprod{\affrt }{\affrt}} \in \{0,1 \}.
\] 

Let $\alpha \in \Phi$ and $k\in \Z$.  Let $H_{\alpha ,k}\subset V$ be the
affine hyperplane
\begin{equation}\label{eqn:affhyperplane}
H_{\alpha ,k} = \{ v\in V\mid \innprod{\alpha}{v} = k\},
\end{equation}
and let $H^{1}_{\alpha ,k}$ be the subset 
\begin{equation}\label{eqn:hone}
H^{1}_{\alpha ,k} = \{v\in V\mid k\leq \innprod{\alpha}{v}\leq k+1 \}.
\end{equation}
Attached to the hyperplane $H_{\alpha , k}$ is the affine reflection
$s_{\alpha ,k}$, which acts on $V$ by
\[
s_{\alpha ,k}\colon v\longmapsto v-2 (\innprod{\alpha}{v}-k)\alpha /\innprod{\alpha}{\alpha}. 
\]
The affine Weyl group $\Waff$ corresponding to $\Phi $ is the group
generated by all the affine reflections $s_{\alpha ,k}$.  We can
represent $\Waff$ as a finitely generated Coxeter group by using the
generating set $\Saff = S \cup \{s_{\affrt , 1} \}$, where we identify
$s_{\alpha}\in S$ with $s_{\alpha ,0}$.

Let $\sH$ be the affine hyperplane arrangement consisting of all
affine hyperplanes of the form \eqref{eqn:affhyperplane}.  The
connected components of $V\smallsetminus \sH$
are called \emph{alcoves}.  There is a distinguished alcove $A_{0}$
defined by 
\[
A_{0} = \{v\mid 0<\innprod{\alpha}{v}<1\quad \text{for all $\alpha \in
\Phi^{+}$}\}, 
\]
and $w\mapsto A_{0}\cdot w$ gives a bijection between $\Waff$ and the
set of alcoves.  We often identify alcoves and elements of $\Waff$
under this bijection.

Any $w\in \Waff$ determines a function $b_{w}\colon
\Phi^{+}\rightarrow \Z$ as follows.  The closure of the alcove
$A_{0}\cdot w$ can be uniquely written as the intersection of subsets
of the form \eqref{eqn:hone}:
\[
\overline{A_0\cdot w} = \bigcap_{\alpha \in \Phi^{+}} H^{1}_{\alpha ,
k (\alpha)}.
\]
We put $b_{w} (\alpha) = k (\alpha) \in \Z$.

An \emph{expression} for $w\in \Waff$ is a representation of $w$
as a product of elements of $\Saff$.  An expression is \emph{reduced}
if it has minimal length among all expressions for $w$.  We define
the length $\ell (w)$ of $w$ to be the length of a reduced expression
for $w$.  

Given an expression $s_{1}\dotsb s_{N}$, a \emph{subexpression} is a
(possibly empty) expression of the form $s_{i_{1}}\dotsb s_{i_{M}}$,
where $1\leq i_{1}<\dotsb <i_{M}\leq N$.  We endow $\Waff$ with a
partial order by defining $u\leq w$ if an expression for $u$ appears
as a subexpression of a reduced expression for $w$.

\subsection{}\label{ss:celldefs}
Next we describe Kazhdan--Lusztig
cells.  Given $w\in \Waff$, we define the \emph{left descent set} $L
(w) \subset \Saff $ by
\[
L (w) = \{ s\in \Saff \mid \ell (sw) < \ell (w) \},
\]
and analogously define the \emph{right descent set} $R (w)$ by the condition
$\ell (ws)<\ell (w)$. 

Let $u,w\in \Waff $, and let $P_{u,w} (t)\in \Z [t]$ be the
\emph{Kazhdan--Lusztig polynomial} attached to the pair $(u,w)$ \cite{kl}.  We
do not recall the definition here, but only mention the following
properties:
\begin{enumerate}
\item $P_{u,w} = 0$
unless $u\leq w$,
\item $P_{u,u} = 1$, and 
\item $P_{u,w}$ has degree at most $d
(u,w) := (\ell (w)-\ell (u)-1)/2$.
\end{enumerate}
We write $\join{u}{w}$ if $u<w$ and $\deg P_{u,w} = d (u,w)$.  If $w< u$ we write
$\join{u}{w}$ if $\join{w}{u}$ holds.

We are now ready to define the cells of $\Waff $.  The \emph{left
$W$-graph} $\Gamma_{L}$ of $\Waff $ is the directed graph with vertex
set $\Waff $, and with an arrow from $u$ to $w$ if and only if
$\join{u}{w}$ and $L (u)\not\subset L (w)$.  We can similarly define
the \emph{right $W$-graph} $\Gamma_{R}$.  Then the left and right
cells of $\Waff$ are extracted from the graphs $\Gamma_{L},
\Gamma_{R}$ as follows.  Recall that given any directed graph, we say
two vertices are in the same \emph{strong connected component} if
there exist directed paths from each vertex to the other.

\begin{definition}\label{def:cells}
The \emph{left cells} (respectively, \emph{right cells}) of $\Waff $
are the strong connected components of the graph $\Gamma_{L}$ (resp.,
$\Gamma_{R}$).  The elements $u,w\in \Waff$ are in the same
\emph{two-sided cell} if they are in the same left \emph{or} right
cell.
\end{definition}

\begin{figure}[htb]
\centering
\subfigure[$\ell =1$\label{fig:g2-1}]{\includegraphics[scale=0.6]{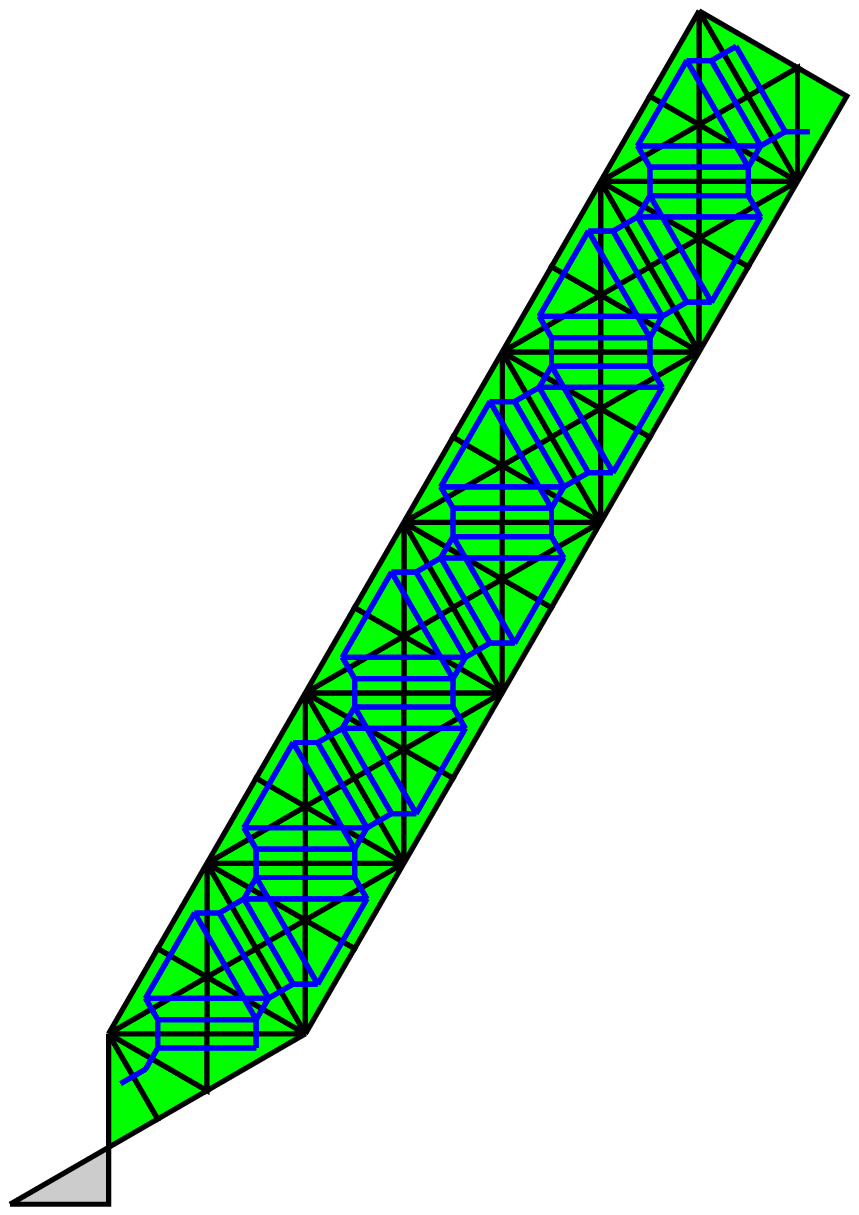}}\quad\quad 
\subfigure[$\ell =3$\label{fig:g2-3}]{\includegraphics[scale=0.6]{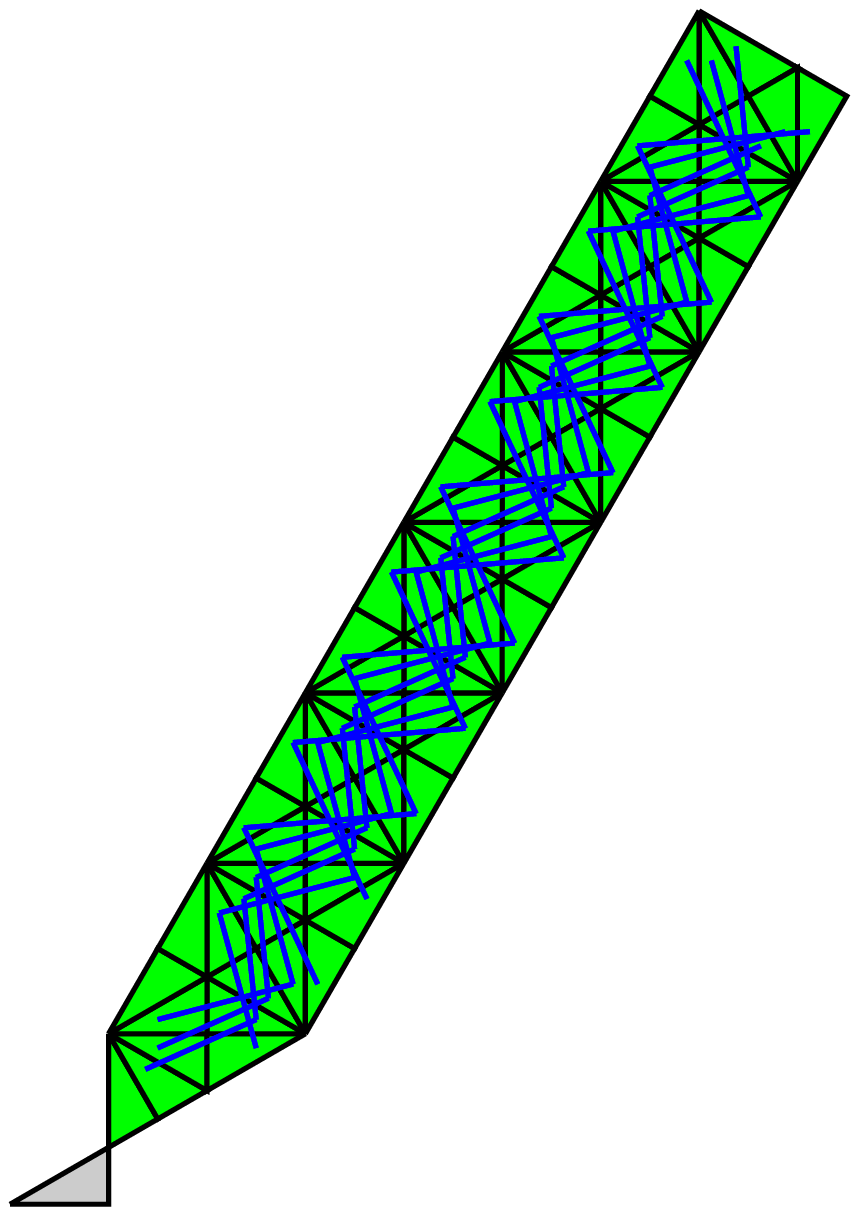}}
\caption{\label{fig:one}}
\end{figure}

\begin{figure}[htb]
\centering
\subfigure[$\ell =5$\label{fig:g2-5}]{\includegraphics[scale=0.6]{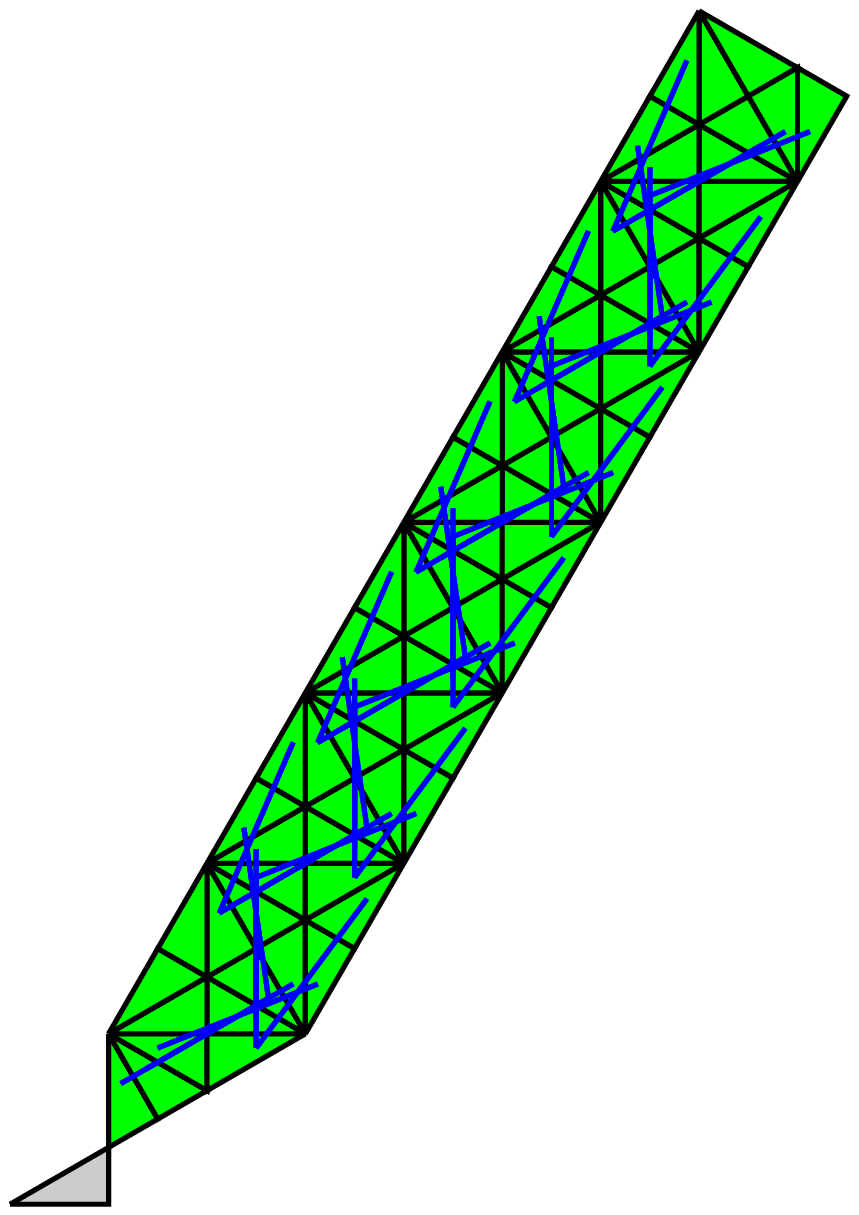}}\quad \quad 
\subfigure[$\ell =7$\label{fig:g2-7}]{\includegraphics[scale=0.6]{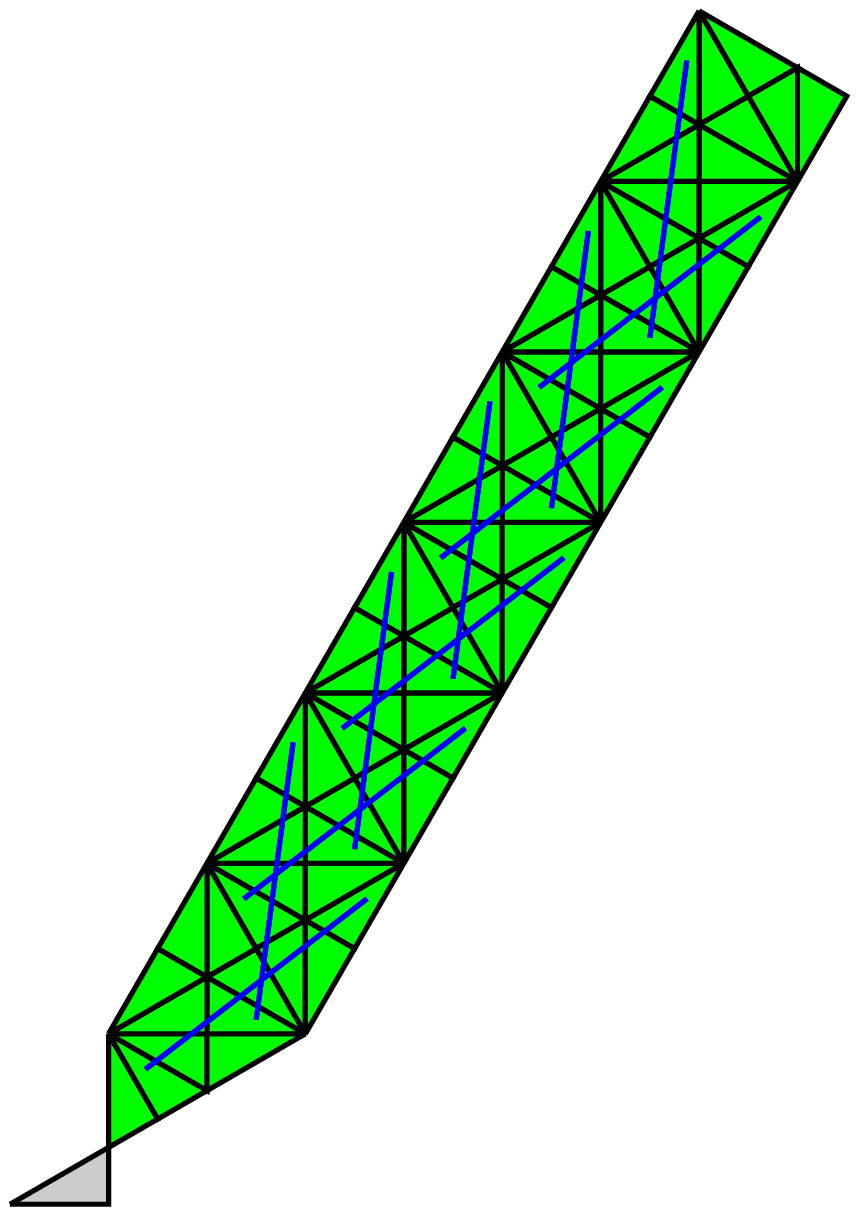}}
\caption{\label{fig:two}}
\end{figure}

\begin{figure}[htb]
\centering
\subfigure[$\ell =9$\label{fig:g2-9}]{\includegraphics[scale=0.6]{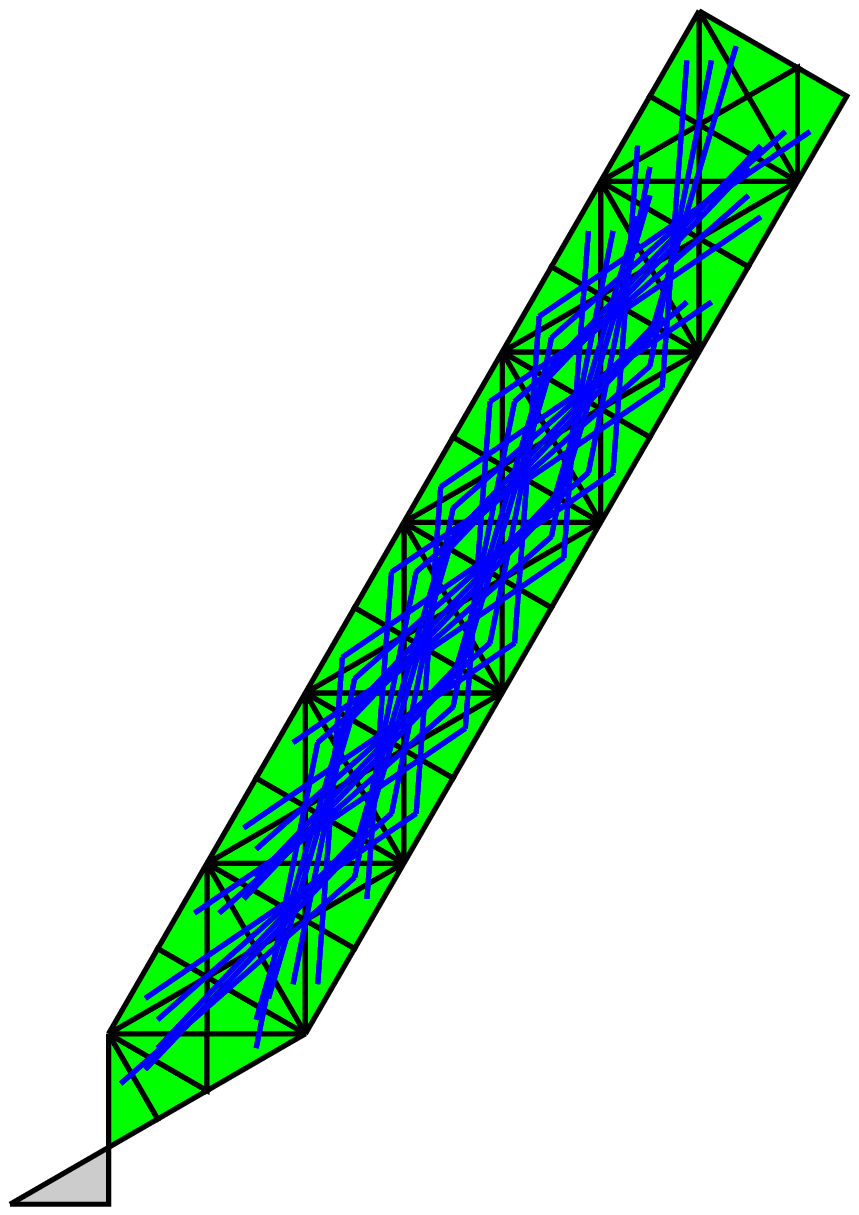}}\quad\quad 
\subfigure[$\ell =11$\label{fig:g2-11}]{\includegraphics[scale=0.6]{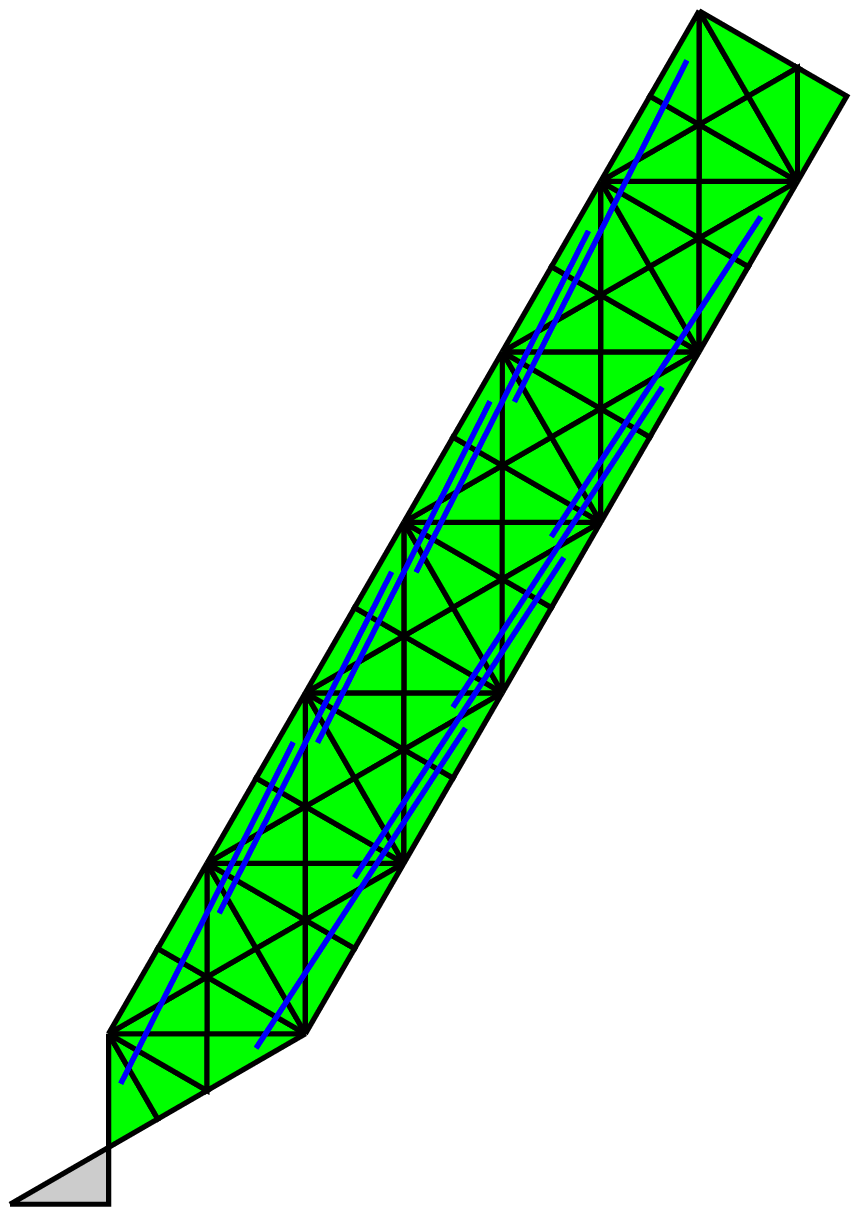}}
\caption{\label{fig:three}}
\end{figure}

It is known that each affine Weyl group has only finitely many
two-sided cells, and that each two-sided cell is a union of finitely
many left cells \cite{cellsII}.  We remark that in general the graphs
$\Gamma_{L}$, $\Gamma_{R}$ are extremely complicated.  Figures
\ref{fig:one}--\ref{fig:three} show the subgraph of $\Gamma_{L}$
corresponding to a particular left cell of $\tilde{G}_{2}$.  The labels
of the figure show the length difference $\ell = \ell (w) - \ell (u)$ of the words
connected by edges in the graph; we have omitted the arrowheads for
clarity.

\subsection{}\label{ss:automata}
Finally we discuss automata.  Let $A$ be a finite set,
and let $A^{*}$ be the free monoid generated by $A$ under the
operation of concatenation.  Thus $A^{*}$ consists of all sequences
$a_{1}\dotsb a_{k}$, where $a_{i}\in A$, together with the empty
sequence.  A \emph{formal language} $L$ is any subset of $A^{*}$.  We
call $A$ the \emph{alphabet} for the language $L$.  Among all formal
languages, our interest lies in the regular languages, which are
defined using finite state automata:

\begin{definition}\label{def:fsa}
A \emph{finite state automaton} $\sA$ over the alphabet $A$ consists of the
following data:
\begin{enumerate}
\item a finite set $Q$, called the set of \emph{states},
\item a unique \emph{initial state} $q_{0}\in Q$,
\item a subset $F\subset Q$, called the set of \emph{accepting states},
and 
\item a function $t\colon Q\times A \rightarrow Q\cup \{\emptyset \}$,
called the \emph{transition function}.  
\end{enumerate}

\end{definition}

A word $a_{1}\dotsb a_{k}\in A^{*}$ is said to be \emph{accepted} by $\sA$ if
there exist states $q_{0},\dotsc,q_{k}$ such that $t (q_{i}, a_{i+1})
= q_{i+1}$ for $0\leq i\leq k-1$, and $q_{k}\in F$.  A language $L$ is
\emph{regular} if there exists a finite state automaton accepting exactly
the words in $L$.

We may picture $\sA$ as a decorated directed graph as follows.  We
label the vertices of the graph by $Q$, and circle the vertices
corresponding to the final states.  We put an edge from $q$ to $q'$
labelled by $a\in A$ if $t (q,a) = q'$.  A path in this graph starting
at $q_{0}$ determines a word in $A^{*}$ by concatenating the edge
labels as one traverses the path.  A word is accepted by $\sA$ if and
only if one can find a path building the word that ends in an
accepting state.

%%%%%%%%%%%%%%
%            %
% automaton  %
%            %
%%%%%%%%%%%%%%

\section{The automaton}\label{s:automaton}

\subsection{}
Let $\Reduced (\Waff)$ be the set of all reduced expressions for all
elements $w\in \Waff$, regarded as a language on the alphabet $\Saff$.
More generally, if $U$ is any subset of $\Waff$, we let $\Reduced (U)$
be the set of all reduced expressions for all $w\in U$. 

By work of Brink--Howlett \cite{bh}, it is known that $\Reduced
(\Waff)$ is regular.\footnote{In fact \cite{bh} shows that the
language of reduced expressions is regular for any Coxeter group.}
Our goal is to show the regularity of the language $\Reduced (C)$ for
any cell of $\Waff$.  To do this we generalize an automaton for
affine Weyl groups first described by Eriksson \cite{eriksson} and
Headley \cite{headley}.

\subsection{} Let $\sH\subset V$ be a finite affine hyperplane
arrangement, and let $\sR$ be the set of connected components of
$V\smallsetminus \sH$.  We will always assume that $\sH$ contains the
hyperplanes $H_{\affrt , 1}$ and $H_{\alpha ,0}$, $\alpha \in \Delta$,
which we denote by $\{H_{s}\mid s\in \Saff \}$.  We say the set of
regions $\sR$ has \emph{property $(*)$} if it satisfies the following
condition:
\begin{quote}
If $R\in \sR$ and $A_{0}$ lie on the same side of $H_{s}$, then there
is a unique $R'\in \sR$ such that $R\cdot s \subset R'$.
\end{quote}

The following result explains the connection between property $(*)$
and automata:

\begin{proposition}\label{prop:starautomata}
Suppose the set of connected components $\sR$ of a finite affine
hyperplane arrangement $\sH$ satisfies $(*)$.  Then there is a finite
state automaton $\sA$ accepting $\Reduced (\Waff)$ with states given
by $\sR$.
\end{proposition}

\begin{proof}
We define $\sA$ by taking its states to be $\sR$ and its initial state
to be $A_{0}$; note that $A_{0}\in \sR$ since $\sH$ contains
$\{H_{s}\mid s\in \Saff \}$.  We declare all states to be accepting.
We define the transition function $t$ by $t (R,s) = \emptyset$ unless
$A_{0}$ and $R$ lie on the same side of $H_{s}$, in which case we put
$t (R,s) = R'$, where $R'$ is the unique region with $R\cdot s \subset
R'$.

The proof that $\sA$ accepts $\Reduced (\Waff)$ is essentially the
same as that for Theorem V.6 in \cite{headley}.  For the convenience of
the reader, we give the details.  First, it is clear that the unique
expression for the identity in $\Waff$ is accepted.  Also all
expressions of length $1$ are accepted, and such expressions are
automatically reduced.

Now suppose that the reduced expression $w=s_{k}\dotsb s_{1}$ is
accepted and let $s\in \Saff$.  We must show that $s_{k}\dotsb s_{1}s$
is accepted if and only if it is reduced.  Recall that $s_{k}\dotsb
s_{1}$ is reduced if and only if the hyperplanes
\begin{equation}\label{eqn:hyplist}
H_{1}, H_{2}\cdot s_{1}, H_{3}\cdot s_{2}s_{1},\dotsc ,H_{k}\cdot
s_{k-1}\dotsb s_{1}
\end{equation}
are distinct and separate $A_{0}$ from $A_{0}\cdot w$
(cf. \cite[\S4.5]{humph.book}; here we abbreviate $H_{s_{i}}$ by
$H_{i}$).  Thus if $ws$ is reduced, then $H_{s}$ must separate $A_{0}$
from $A_{0}\cdot ws$.  This means $H_{s}$ does not separate $A_{0}$
from $A_{0}\cdot w$, which implies $s_{k}\dotsb s_{1}s$ is accepted by
$\sA$.  

Conversely, suppose $s_{k}\dotsb s_{1}s$ is accepted by $\sA$.  Then
$H_{s}$ does not separate $A_{0}$ from $A_{0}\cdot w$, which means
$H_{s}$ separates $A_{0}$ from $A_{0}\cdot ws$.  Certainly the
hyperplanes obtained by applying $s$ (on the right) to
\eqref{eqn:hyplist} are distinct, and none equal $H_{s}$.  Moreover if
some hyperplane $H = H_{j}\cdot s_{j-1}\dotsb s_{1}s$ from this new
list fails to separate $A_{0}$ from $A_{0}\cdot ws$, then $H\cdot s$
fails to separate $A_{0}$ from $A_{0}\cdot ws$, and thus cannot be
part of \eqref{eqn:hyplist}.  This completes the proof.
\end{proof}

The automaton from Proposition \ref{prop:starautomata} also satisfies
the following property: if an expression $w=s_{1}\dotsb s_{k}$ is
reduced and is accepted by the state $R$, then the alcove $A_{0}\cdot
w$ lies in the region $R$.

\subsection{}
Now let $N\geq 0$ be an integer, and let $\sH_N$ be the arrangement
\[
\sH_N = \{H_{\alpha ,k}\mid \alpha \in \Phi^{+}, k = -N,\dotsb N+1 \}.
\]
Let $\sR_{N} = V\smallsetminus \sH_{N}$.  Then we have the following theorem:

\begin{theorem}\label{thm:affineaut}
The set of regions $\sR_{N}$ has the property $(*)$, and thus the
automaton built from $\sR_{N}$ using Proposition
\ref{prop:starautomata} accepts the language $\Reduced (\Waff)$.
\end{theorem}

\begin{proof}
The proof is a generalization of \cite[Lemma V.5]{headley}.  For each
$\alpha \in \Phi^{+}$, and for $k=-N-1,\dotsc ,N+1$ define ``root
strips'' by
\[
R^k_{\alpha } = \begin{cases}
\{v\mid N +1< \innprod{\alpha}{v} \}&k=N+1,\\
\{v\mid \innprod{\alpha}{v} <-N \}&k=-N-1,\\
\{v\mid k<\innprod{\alpha}{v}<k+1 \}&\text{otherwise}.
\end{cases}
\]
The elements of $\sR_{N}$ are just the connected components of
all possible intersections of the $R^{k}_{\alpha}$.

Let $R\in \sR_{N}$, and suppose $R$ and $A_{0}$ lie on the same side
of $H_{s}$.  We must show that $R\cdot s$ is contained in a unique element
of $\sR_{N}$.  By the above discussion, we know that for each $\alpha
\in \Phi^{+}$, we have $R\subset R^{k}_{\alpha}$ for some $k$
depending on $\alpha$.  We must show that $R\cdot s \subset R^{l}_{\beta}$
for some $\beta \in \Phi^{+}$ and some $l$.

First assume $s$ is not the affine reflection $s_{\affrt ,1}$.  Then
either $\alpha =\alpha_{s}$ or $\alpha \not = \alpha_{s}$.  In the
first case, if $R\subset R^{k}_{\alpha}$ then we must have $k\geq 0$,
and it is clear that $R\cdot s \subset R_{\alpha}^{-k-1}$.  Note that
the hypothesis that $R$ and $A_{0}$ lie on the same side of $H_{s}$ is
essential, since $R_{\alpha}^{-N-1}\cdot s$ meets both
$R_{\alpha}^{N}$ and $R_{\alpha}^{N+1}$.

Next suppose $\alpha \not =\alpha_{s}$.
Then if $k\not = N+1, -N-1$, we have
\begin{align*}
R^{k}_{\alpha}\cdot s &= \{v\mid k <\innprod{\alpha\cdot s}{v}<k+1 \} \\
	&=R^{k}_{\alpha\cdot s}.
\end{align*}
Moreover, $\alpha\cdot s\in \Phi^{+}\smallsetminus \{\alpha_{s} \}$,
since the only positive root $s$ makes negative is $\alpha_{s}$.
Hence in this case $R\cdot s$ is taken into a unique region of $\sR_{N}$ if
$s$ is not the affine reflection; if $k=N+1$ or $-N-1$ we argue
similarly.

Finally assume $s$ is the affine reflection, and suppose $R\subset
R^{k}_{\alpha}$.  Again there are two possibilities to consider.  If
$\alpha$ is the highest root $\affrt$, then since $R$ lies on the same
side of $H_{s}$ as $A_{0}$, we must have $k\leq 0$.  Therefore
$R\cdot s\subset R_{\alpha}^{-k+1}$.  On the other hand, if $\alpha$ is not
the highest root, then $c =
2{\innprod{\alpha}{\affrt}}/{\innprod{\affrt}{\affrt}}$ either equals
$0$ or $1$.  If $c=0$, then $s_{\affrt ,1}$ stabilizes all the sets
$R_{\alpha}^{k}$.  If $c=1$, then $s_{\affrt ,1} (\alpha) = \alpha
-\affrt$, which is negative.  We have for $k\not = N+1$,
\begin{align*}
R^{k}_{\alpha}\cdot s &= \{v\mid k< \innprod{\alpha -\affrt}{v} <k+1 \} \\
		&= \{v\mid -k-1 < \innprod{\alpha '}{v}<-k \}\\
		&= R_{\alpha '}^{-k-1},
\end{align*}
where $\alpha '\in \weylposrts$ is the root $\affrt-\alpha $.  If
$k=N+1$, then $R^{N+1}_{\alpha}\cdot s \subset R_{\alpha '}^{-N-1}$.
This completes the proof.
\end{proof}

\begin{remark}\label{rem:modified}
Recall that $\Phi$ is \emph{simply-laced} if all roots have the same
length.  Suppose $\Phi$ is not simply-laced, and partition the
positive roots $\Phi^{+}$ into the long and short roots
$\Phi^{+}_{l}\cup \Phi^{+}_{s}$.  Choose a function $\nu \colon
\Phi^{+}\rightarrow \Z_{\geq 0}$ that is constant on $\Phi^{+}_{l}$,
$\Phi_{s}^{+}$, and consider the affine arrangement $\sH_{\nu}$
consisting of the hyperplanes
\begin{equation}\label{eqn:hyp}
\sH_{\nu} = \{H_{\alpha ,k} \mid k = -\nu (\alpha),\dotsc ,\nu
(\alpha)+1 \}.
\end{equation}
Then the proof of Theorem \ref{thm:affineaut} actually shows that the
set of regions of the complement $V\smallsetminus \sH_{\nu}$ has
property $(*)$ as well, and thus can be used as the set of states for
an automaton recognizing $\Reduced (\tilde{W})$.  Indeed, the only
observations one needs to make are that the action of $W$ preserves
root lengths, and that if $2
\innprod{\alpha}{\affrt}/\innprod{\affrt}{\affrt} = 1$ then $\alpha
-\affrt$ has the same length as $\alpha$.
\end{remark}

\begin{figure}[htb]
\begin{center}
\includegraphics[scale=0.7]{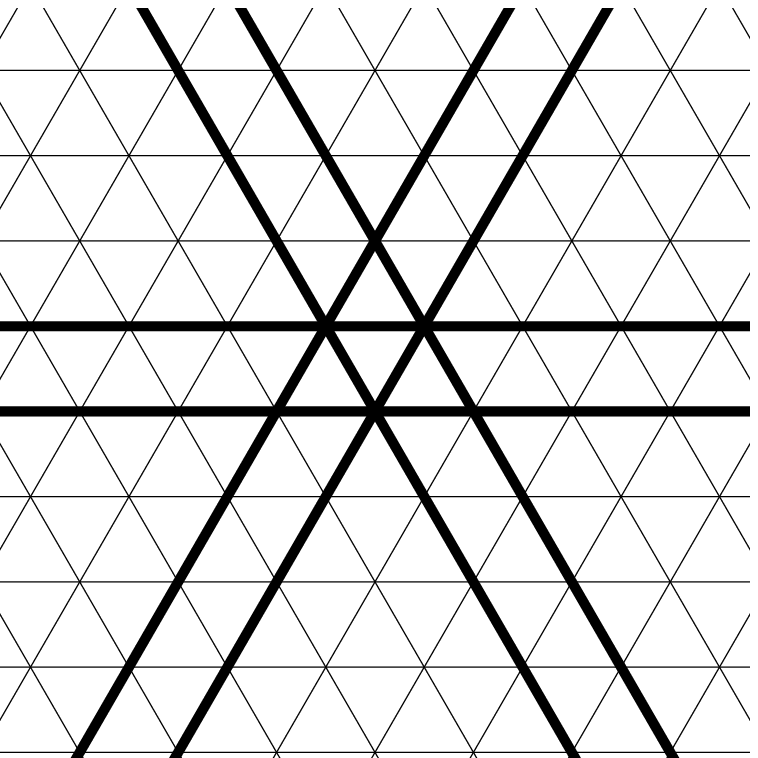}
\end{center}
\caption{\label{fig:A2shi}The regions in $\sR_{0}$ for $\tilde{A}_{2}$}
\end{figure}

\begin{remark}
The arrangement $\sH_{N}$ is known in the combinatorics literature as
the \emph{extended Shi arrangement} \cite{stanley}.  The case $N=0$
was studied by Shi \cite{shi.sign}, who called them \emph{sign-type
regions}.  Figure \ref{fig:A2shi} shows these regions and the alcoves
for $\tilde{A}_{2}$.  Shi showed that the number of regions in
$\sR_{0}$ equals $(h+1)^{r}$, where $h$ is the Coxeter number of
$\Phi$ and $r$ is the rank of $\Phi$, and also used sign-type regions
to explicitly describe the Kazhdan--Lusztig cells for $\tilde{A}_{n}$
\cite{shi}.  For the extended Shi arrangement, Athanasiadis \cite{ath}
showed that the number of regions in $\sR_{N-1}$ is $(Nh+1)^{r}$ if
$\Phi$ is of classical type, i.e.~type $A$, $B$, $C$, $D$.
\end{remark}

%%%%%%%%%%%%%%%%%%%%%%%%%%%
%                         %
%  Kazhdan-Lusztig cells  %
%                         %
%%%%%%%%%%%%%%%%%%%%%%%%%%%

\section{Kazhdan-Lusztig cells and automata}\label{s:cells}

\subsection{}
In this section we prove Theorem \ref{thm:main}.  The main idea of the
proof is to show that we can choose $N$ sufficiently large such that
for any left cell $C$, we can find a finite set of regions
$\{R_{i}\}\subset \sR_{N}$ such that the alcoves in the union of the
$R_{i}$ are exactly those in $C$.  The automaton accepting $\Reduced
(C)$ is then given by modifying the automaton from Theorem
\ref{thm:affineaut} to make only the states corresponding to the
$R_{i}$ accepting.

\begin{figure}[htb]
\begin{center}
\includegraphics[scale=0.8]{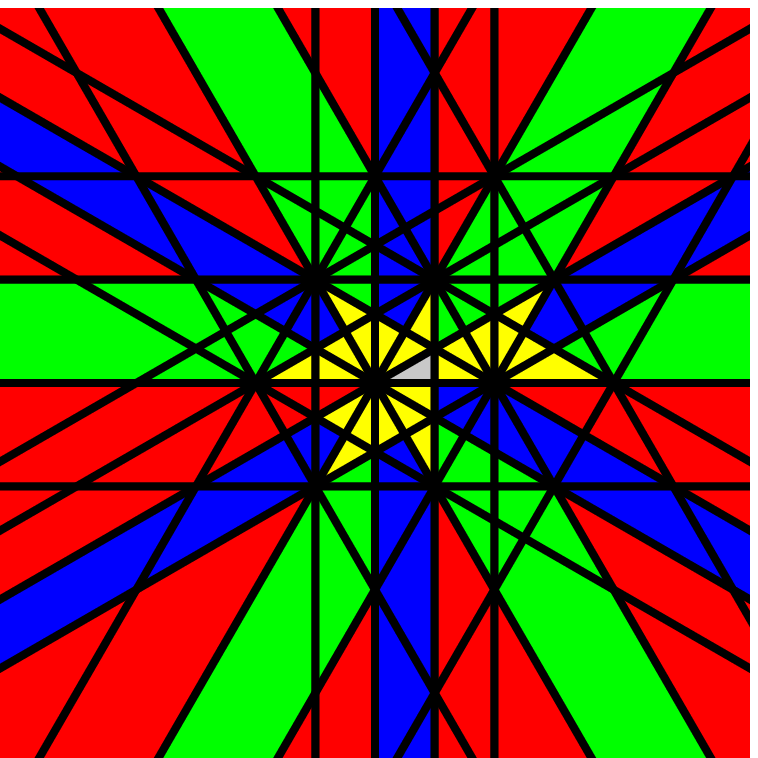}
\end{center}
\caption{\label{fig:g2}The left cells of $\tilde{G}_{2}$ and the regions in
$\sR_{1}$}
\end{figure}

Figure \ref{fig:g2} shows an example for $\Waff = \tilde{G}_{2}$.  The
dark lines are the hyperplanes in $\sH_{1}$; each region of a given
color is a two-sided cell, and the connected regions of a given color
are the left cells.  Note that we have not drawn the alcoves.
It is clear from the figure that any left cell is a union of
regions from $\sR_{1}$.  Note that $N=1$ is the smallest value we can
take for this to work; in particular left cells in $\tilde{G}_{2}$ are
not unions of sign-type regions.

\subsection{}
The main tool we need to carry out the proof is a result of Du \cite{du.polyhedra}, who
showed that left cells are finite unions of certain polyhedra in $V$.
More precisely, let $\sP$ be the set of all polyhedra in $V$ of full
dimension bounded by finitely many affine hyperplanes of the form
$H_{\alpha , k}$, $\alpha \in \Phi^{+}$, $k\in \Z$.  Then we have the
following theorem:

\begin{theorem}\label{thm:du}
\emph{\cite{du.polyhedra}} Let $C$ be a left cell of $\Waff$.  Identify $C$ with the closure of
the set of alcoves $\{A_{0}\cdot w \mid w\in C\}$ in $V$.  Then
there
exists a finite subset $\sP (C)\subset \sP$ such that 
\[
C = \bigcup_{P\in \sP (C)} P.
\]
\end{theorem}

We give some indications of the proof of Theorem \ref{thm:du}.  Let
$P\in \sP$.  Assume there exists a Weyl chamber $\sC$ containing 
$P$.  Let $\Delta '\subset \Phi$ be the simple system determined by
$\sC$ (so that all $\alpha \in \Delta '$ are nonnegative on $\sC$).
Then the dimension $\dim P$ of $P$ is defined to be the cardinality of
the set\footnote{We remark that there is a typo in
\cite{du.polyhedra}, in which ``unbounded'' is replaced by ``bounded''
in the definition.}
\[
\{\alpha \in \Delta '\mid \text{$\alpha$ is unbounded on $P$} \}.
\] 

Du defines a certain class of polyhedra in $\sP$ called \emph{special
polyhedra} as follows.  Recall \eqref{eqn:hone} that we have defined
the subsets $H^{1}_{\alpha ,k}$.  Fixing $\sC$, we let $H^{+}_{\alpha
,k}$ be the closure of the unique component of $V\smallsetminus
H_{\alpha ,k}$ that meets any translate of $\sC$.  Given $\Delta '$ as
above we choose a function $b\colon \Delta '\rightarrow \Z_{\geq 0}$
and a subset $\Lambda \subset \Delta '$.  Let $\overline{\Lambda } =
\Delta '\smallsetminus \Lambda$.  Then the special polyhedron $P
(\Lambda ,b)$ is defined by
\[
P (\Lambda ,b) = \bigl(\bigcap_{\alpha \in \Lambda} H^{1}_{\alpha ,b
(\alpha)}\bigr) \cap \bigl(\bigcap_{\alpha \in \overline{\Lambda}}
H^{+}_{\alpha ,b (\alpha)}\bigr).
\]
Clearly any special polyhedron lies in $\sP$, and so does any
polyhedron built by forming intersections of special polyhedra with
subsets of the form $H^{1}_{\alpha,k}$ and $H^{+}_{\alpha,k}$.  Du
proves the following properties of special polyhedra: 

\begin{enumerate}
\item For any special polyhedron $P$, we can find a subpolyhedron
$P'\subset P$, also special, with $P\smallsetminus P'$ a finite union
of special polyhedra $P_{i}$ of lower dimension than that of $P$.
\item The polyhedron $P'$ can be chosen such that $P'$ is a finite
union of polyhedra $Q_{j}\in \sP$ (not special in general), and such
that each $Q_{j}$ is contained in a single left cell.  These $Q_{j}$
are built using intersections of $P$ with subsets of the form
$H^{1}_{\alpha,k}$ and $H^{+}_{\alpha,k}$.
\end{enumerate}

Finally Du completes the proof by induction and by choosing a finite
collection of special polyhedra such that all alcoves will eventually
be accounted for by the $Q_{j}$ constructed in property (2).  For this
he uses the sign-type regions $\sR_{0}$.  Given any sign-type region
$X$, he shows that there is a canonical special polyhedron $P (X)
\supset X$ attached to $X$ with $\dim X = \dim P (X)$, and such that
$P (X)\smallsetminus X$ is a finite union of sign-type regions.  The
sign-type regions account for all alcoves, and thus so do the special
polyhedra $P (X)$.  Using induction on $\dim X$ and properties (1),(2)
completes the proof.

The most subtle part of the proof of Theorem \ref{thm:du} is, given a
special polyhedron $P$, the construction of $P'$ and the polyhedra
$Q_{j}$.  For this Du constructs certain infinite sequences $y_{1},
y_{2}, \dotsc, y_{i}, \dotsc $ in special polyhedra that become
left-equivalent for $i$ sufficiently large.  A key ingredient here is
the boundedness of Lusztig's $a$-function for affine Weyl groups
\cite{lusztig}.

Returning to the general discussion, here is the connection between
polyhedra from $\sP$ and the regions $\sR_{N}$:

\begin{lemma}\label{lem:bigN}
Let $\{P_{i} \}$ be a finite subset of $\sP$, and let $P$ be the (not
necessarily convex) union of the $P_{i}$.  Then for $N$ sufficiently
large, we can find a finite set of regions $\sR (P)\subset \sR_{N}$
such that
\[
P = \bigcup_{R\in \sR (P)} \overline{R}.
\]
\end{lemma}

\begin{proof}
First we assume that $P$ is a single polyhedron from $\sP$.  If $P$ is
bounded, then $P$ is the closure of a finite union of alcoves.  Since any
finite collection of alcoves can appear as regions in $\sR_{N}$ for
$N\ggt 0$, the statement follows.  

If $P$ is unbounded, then there are only finitely many positive roots
$\alpha$ such that $\innprod{\alpha}{v}\geq k (\alpha)$ is a defining
inequality for $P$, where $k(\alpha) \in \Z$.  Then we simply choose
\[
N\ggt \max_{\alpha \in \Phi^{+}} |k (\alpha)|,
\]
and $P$ can be written as a union of
regions from $\sR_{N}$.

Finally, if $P = \cup P_{i}$ is a union of polyhedra, then we compute
$N_{i}$ for each $P_{i}$ as above and then choose $N\ggt \max N_{i}$.
This completes the proof.
\end{proof}

\subsection{}
We are now ready to prove our main theorem:

\begin{theorem}\label{thm:main}
Let $C$ be a left cell in $\Waff$, and let $\Reduced (C)$ be the
language of all reduced expressions of $w\in C$.  Then $\Reduced (C)$
is regular.
\end{theorem}

\begin{proof}
By Theorem \ref{thm:du}, the cell $C$ is a union $P=\cup P_{i}$ of
finitely many polyhedra from $\sP$.  By Lemma \ref{lem:bigN}, we may
find $N\ggt 0$ such that each $P_{i}$ is a union of the closures of
regions of $\sR_{N}$.  Let $\sA_{N}$ be the automaton from Theorem
\ref{thm:affineaut} constructed from $\sR_{N}$.  We define a finite
state automaton $\sA (C)$ by modifying $\sA_{N}$ so that the only
accepting states are those corresponding to the regions of $\sR_{N}$
in $P$.  By the comment after the proof of Proposition
\ref{prop:starautomata}, an expression $w=s_{1}\dotsb s_{k}$ is
accepted by $\sA (C)$ if and only if it is reduced and $w\in C$.  This
completes the proof.
\end{proof}

As remarked before, each two-sided cell contains only finitely many
left cells \cite{cellsII}.  Moreover, any right cell $C_{R}$ is
obtained by inverting some left cell $C$, which means $\Reduced
(C_{R})$ consists of all words in $\Reduced (C)$ reversed.  Since both 
a finite union of regular languages and the reversal of a regular
language are regular, we obtain the following corollary:

\begin{corollary}\label{cor:allcells}
Let $C$ be any cell in $\Waff$, including right and two-sided.  Then
$\Reduced (C)$ is regular.
\end{corollary}

\begin{remark}
It is interesting to consider the minimal value of $N$ needed to
simultaneously show that all the cells of $\Waff$ give regular
languages.  A consequence of Shi's work \cite{shi} is that $N=0$
suffices for type $A$, and Figure \ref{fig:g2} shows that $N=1$ is the
minimal value needed for $\tilde{G}_{2}$.  Figure \ref{fig:g2} also
shows that the arrangement $\sH_{1}$ contains extra lines not needed
to decompose left cells into regions.  For instance, the left- and
rightmost vertical lines, as well as the highest and lowest horizontal
lines, are not supporting boundary lines for any left cell, and thus
are not needed to distinguish left cells from each other.    

Thus one has the natural question of defining the smallest possible
automaton needed to describe all the left cells of $\Waff$.  Such an
automaton would have theoretical value, since using it one could
attempt to extend Shi's work for type $A$ \cite{shi} to all types,
especially the classical types (cf. \cite{bedard1, du, du2, chen}).
The obvious idea is to eliminate some hyperplanes from $\sH_{N}$, such
as the aforementioned lines in Figure \ref{fig:g2}.  However, one must
be careful to ensure that the resulting set of regions still satisfies
property $(*)$.  For instance, in Figure \ref{fig:g2} one could try
eliminating the outermost vertical lines, since they are clearly not
needed to distinguish left cells.  But then to keep property $(*)$ one
will be forced to remove other lines from the arrangement, and some of
these lines are necessary to separate left cells.

Hence one must compromise: some extra hyperplanes can be deleted, but
some must remain to preserve property $(*)$.  Based on Remark
\ref{rem:modified} and examples, we propose the following conjecture:

\begin{conjecture}\label{conj:gensigntypes}
Define $\nu \colon \Phi^{+}\rightarrow \Z_{\geq 0}$ by $\nu
(\alpha)=0$ if $\alpha$ is short and $\nu (\alpha)=1$ if $\alpha$ is
long. Let $\sH_{\nu}$ be the affine hyperplane arrangement from
\eqref{eqn:hyp}, and let $\sR_{\nu}$ be the connected components of the
complement $V\smallsetminus \sH_{\nu}$.  Then any left cell of $\Waff$
is a union of regions from $\sR_{\nu}$.
\end{conjecture}

We have checked Conjecture \ref{conj:gensigntypes} for
$\tilde{C}_{2}$ (Figure \ref{fig:c2}), $\tilde{G}_{2}$ (Figure
\ref{fig:g2modified}), $\tilde{C}_{3}$, and the canonical left cells
of $\tilde{B}_{3}$.

\begin{figure}[htb]
\centering
\subfigure[$\tilde{C}_{2}$\label{fig:c2}]{\includegraphics[scale=0.7]{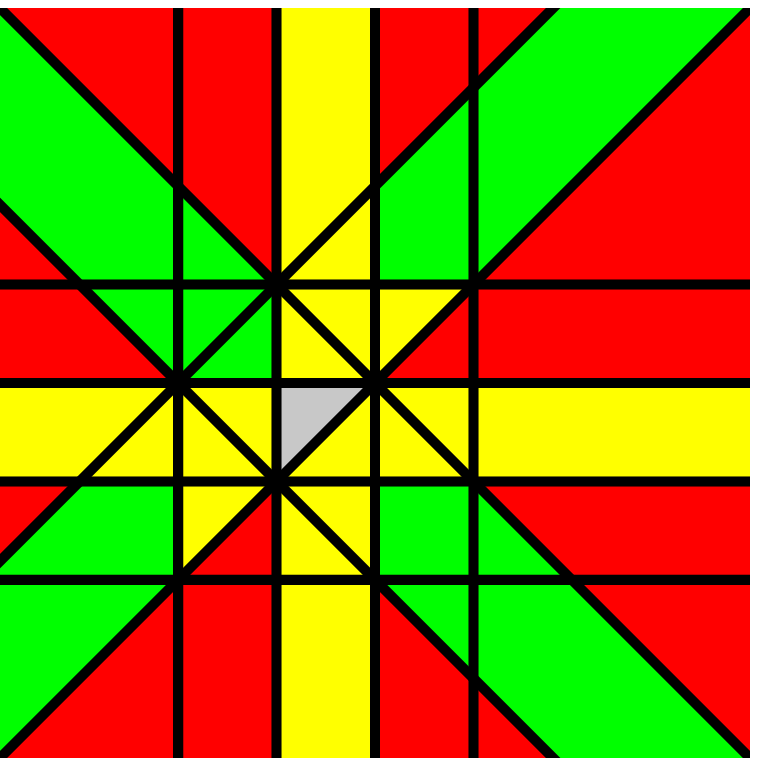}}
\quad\quad
\subfigure[$\tilde{G}_{2}$\label{fig:g2modified}]{\includegraphics[scale=0.7]{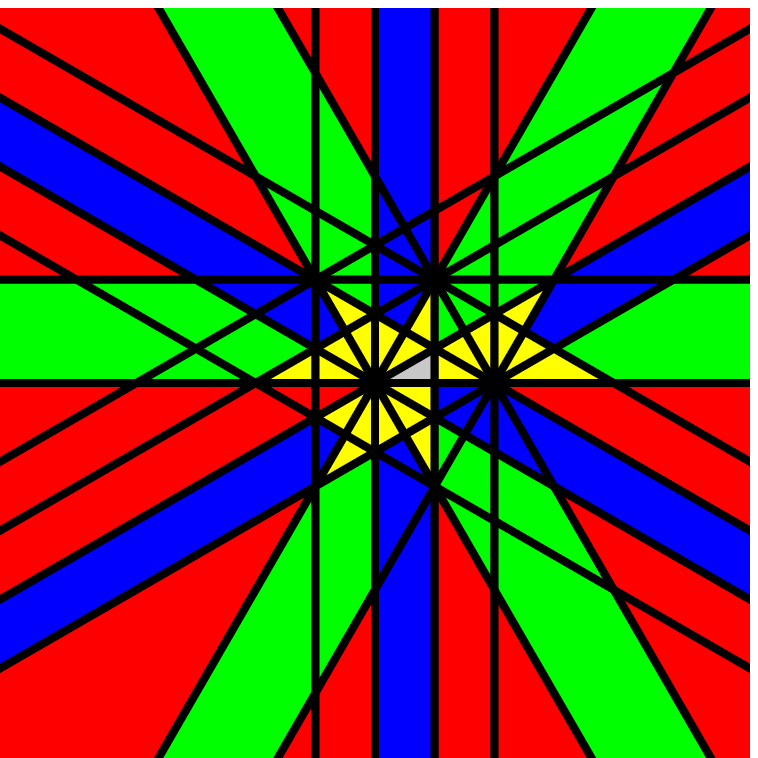}}
\caption{\label{fig2}}
\end{figure}

\end{remark}

\begin{remark}
One can ask if Theorem \ref{thm:main} and the stronger Conjecture
\ref{conj:gensigntypes} apply to the case of cells with \emph{unequal
parameters} \cite{lusztig.unequal}.  There is some evidence that both
are true, although this evidence is limited since there are fewer
cases where one knows all the cells in complete detail.  J.~Guilhot
\cite{guilhot} has recently given a conjectural description of the
left cells in $\tilde{G}_{2}$ for all choices of parameters, and has
proved them in many cases.  For all the examples in his conjectural
description, Conjecture \ref{conj:gensigntypes} holds.  
\end{remark}

\begin{remark}
It is natural to ask if the ideas in this paper can be used to prove
Conjecture \ref{conj:casselman} for all Coxeter groups.  For a general
group, the analogue of the hyperplane arrangement $\sH_{0}$ is the
arrangement of \emph{minimal} or \emph{small} roots.  This is a
certain finite subset of the set of all root hyperplanes in the
geometric realization of $W$ that plays the decisive role in proving
the regularity of $\Reduced (W)$.  Then one could try to identify a
large finite subset of hyperplanes containing the minimal ones that
could be used to determine an analogous class of polyhedra $\sP$.
Unfortunately this naive generaliztion cannot work, since simple
examples show that there exist left cells bounded by infinitely many
root hyperplanes.
\end{remark}

\bibliographystyle{amsplain_initials}
\bibliography{automata}
\end{document}